\documentclass[11pt,amssymb,amstex]{amsart}
\usepackage{graphicx}
\usepackage{epstopdf}

\newtheorem{theorem}{Theorem}[section]

\newtheorem{definition}[theorem]{Definition}
\newtheorem{remark}[theorem]{Remark}
\numberwithin{equation}{section}

\begin{document}

\title{Mendelian and Non-Mendelian Quadratic Operators}

\author{Nasir Ganikhodjaev}
\address{Nasir Ganikhodjaev\\
Department of Computational \& Theoretical Sciences \\
Faculty of Sciences, International Islamic University Malaysia\\
P.O. Box, 141, 25710, Kuantan\\
Pahang, Malaysia} \email{{\tt nasirgani@@hotmail.com}}

\author{Mansoor Saburov}
\address{Mansoor Saburov\\
Department of Computational \& Theoretical Sciences \\
Faculty of Science, International Islamic University Malaysia\\
P.O. Box, 141, 25710, Kuantan\\
Pahang, Malaysia} \email{{\tt msaburov@@gmail.com}}

\author{Uygun Jamilov}
\address{Uygun Jamilov\\
Institute of Mathematics and Information Technologies, Tashkent, Uzbekistan}
\email{{\tt jamilovu@@yandex.ru}}


%

\begin{abstract}
{In this paper, we attempt to provide mathematical models of Mendelian and Non-Mendelian inheritances of the bisexual population system having Fisher's {\textbf{1:1}} principle. In our model, we always assume that distributions of the same phenotype of female and male populations are equal. We study the evolution of a Mendelian trait. As an application of a non-Mendelian inheritance, we construct a quadratic stochastic operator that describes transmission of {\textbf{ABO}} and Rh blood groups.}


\end{abstract}

\maketitle

\section{Introduction}

The study of sex allocation is often considered the most successful
branch of evolutionary biology. Sex allocation is the allocation of
resources to male and female reproductive function
in sexual species (see \cite{CHar}, \cite{West}).

In anthropology and demography, the human sex ratio is the ratio of females to males in the population system. Like most sexual species, the human sex ratio is approximately {\textbf{1:1}}. {\textit{Fisher's principle}} is an evolutionary model that explains why the sex ratio of most sexual species is approximately {\textbf{1:1}}. It was famously outlined by R. Fisher in his book \cite{Fisher}. A. W. F. Edwards has remarked that {\textit{Fisher's principle}} is "probably the most celebrated argument in evolutionary biology" \cite{Edward}. However, W.D. Hamilton had been introduced a model by breaking the assumptions made in Fisher's model in which the population system has an extraordinary sex ratio  (see \cite{Hamil}).  Sex ratios that are {\textbf{1:1}} are hence known as "\textit{Fisherian}", and those that are not {\textbf{1:1}} are "\textit{non-Fisherian}" or "\textit{extraordinary}".

Heredity is the passing of traits to offspring from its parent or ancestors. This is the process by which an offspring cell or organism acquires or becomes predisposed to the characteristics of its parent cell or organism. Inherited traits are controlled by genes and the complete set of genes within an organism's genome is called {\textit{its genotype}}. The idea of particulate inheritance of genes can be attributed to G. Mendel (see \cite{Henig}). {\textit{Mendelian inheritance}} is a scientific description of how hereditary characteristics are passed from parent organisms to their offspring. In Mendelian inheritance, each parent contributes one of two possible alleles for a trait. In humans, eye color is an inherited characteristic and an individual might inherit the "brown-eye trait" from one of the parents (see \cite{Sturm}). {\textit{Non-Mendelian inheritance}} is a general term that refers to any pattern of inheritance in which traits do not segregate in accordance with Mendel's laws.
Inheritance of traits in {\textit{fungi, viruses, and bacteria}} are all non-Mendelian. Non-Mendelian inheritance plays a role in several disease processes (see \cite{Van}).

A quadratic stochastic operator (QSO) is a primary source for investigations of dynamical properties of population genetics \cite{Br}. The fascinating applications of quadratic stochastic operators to population genetics were given in \cite{Lyu}. It describes a distribution of a species for next generation if given the distribution of these species for current generation. In the paper \cite{GRMFRU}, it was given a long self-contained exposition of the recent achievements and open problems in the theory of quadratic stochastic operators.

In this paper, we are attempting to give a mathematical model of Mendelian and Non-Mendelian inheritances in the bisexual population systems. Our approach providing the evolution of the bisexual population system is totally different from the one given in \cite{Lyu} and it gives an opportunity to formulate the Mendelian inheritance in the system. We shall study the evolution of Mendelian and Non-Mendelian inheritances of the biosphere having Fisher's {\textbf{1:1}} principle. In our model, we always assume that "\emph{the distributions of the same phenotype of female and male populations are equal.}" We construct a quadratic stochastic operator that describe inheritance of {\textbf{ABO}} and Rh blood groups. This study is the continuation of the papers \cite{GN} and\cite{GNJDMU}.

\section{Preliminaries}

We recall some definitions and notions (see \cite{GN}).

Let $(\Lambda,L)$ be a graph without loops and multiple edges, where $\Lambda$ is the finite set of vertexes and $L$ is the set of edges of the graph and $\{\Lambda_i\}, \  \ i=1,...,n-$ the set of maximal connected subgraphs (connected components) of the graph $(\Lambda,L)$. Furthermore, let $\Phi$ be a finite set, called the set of alleles (in problems of statistical mechanics, $\Phi-$ is called the range of spin).

The function $\sigma:\Lambda \rightarrow \Phi$ is called a genotype (in statistical mechanics, it is called configuration).
Denote by $\Omega$ the set of all genotypes and by $S(\Lambda,\Phi)$ the set of all probability distributions
defined on the set $\Omega$. A QSO $V:S(\Lambda,\Phi)\rightarrow S(\Lambda,\Phi)$ is
defined as follows: for an arbitrary measure $\lambda \in S(\Lambda,\Phi)$, the measure $\lambda'=V\lambda$ defined
as
\begin{equation} \label{operator0}
\lambda'(\sigma)=\sum\limits_{\sigma',\sigma''\in \Omega} p_{\sigma'\sigma'',\sigma}
\lambda(\sigma')\lambda(\sigma''), \quad \ \forall \ \sigma\in\Omega
\end{equation}
where $p_{\sigma'\sigma'',\sigma}=p_{\sigma''\sigma',\sigma}\geq0$ and
$\sum\limits_{\sigma \in \Omega} p_{\sigma'\sigma'',\sigma}=1 \ \ $ (see \cite{Br}, \cite{Lyu}). Here $p_{\sigma'\sigma'',\sigma}$ is called a heredity coefficient.

In the model \eqref{operator0}, the gender difference of species has not taken into the consideration. However, in many population systems, this difference plays an important role during the evolution  of the population system. Therefore, in this paper, we are going to provide the model of the bisexual population system in which the gender difference has taken into the consideration. It is worth mentioning that our model for  the bisexual population system is totally different from the model provided in \cite{Lyu} and it has an advantage to provide the mathematical model of a Mendelian inheritance.

Let $(\Lambda,L)$ be a graph and $\{\Lambda_i\},$ where $i=0,1,...,n$, be sets of its maximal connected subgraphs such that $\Lambda_0=\{x^0\}.$  Suppose that the population system has only female and male genders $G=\{f,m\}$. Let $\Phi$ be a finite set of alleles. Every species can be characterized by the following genotype
$\sigma: \Lambda\rightarrow G\times\Phi$ such that $\sigma\mid_{\Lambda_0}:\Lambda_0\to G, $ and $\sigma\mid_{\Lambda\setminus\Lambda_0}:\Lambda\setminus\Lambda_0\to \Phi$. We denote the set of all genotypes by $\Omega$. Let $\Omega_f=\{\sigma\in \Omega: \sigma(x^0)=f\}$ and $ \Omega_m=\{\sigma\in \Omega: \sigma(x^0)=m\}$ be sets of all female and male genotypes, respectively. It is clear that $\Omega=\Omega_f\cup\Omega_m$ and $\Omega_f\cap\Omega_m=\emptyset.$ For the two genotypes (parents) $\sigma',\sigma'' \in \Omega,$ we define a set (traits of children) $\Omega (\sigma',\sigma'')$ such that if $\sigma'\mid_{\Lambda_0} = \sigma''\mid_{\Lambda_0}$ then $\Omega (\sigma',\sigma'')=\emptyset$ and if $\sigma'\mid_{\Lambda_0} \neq \sigma''\mid_{\Lambda_0}$ then $\Omega (\sigma',\sigma'')$ is some nonempty subset of $\Omega.$ From the biological point of view, $\Omega (\sigma',\sigma'')$ is a set of all possible genotypes (traits) of children having parents with genotypes $\sigma'$ and $\sigma''$. In this sense, the constrain for the set $\Omega (\sigma',\sigma'')$ means that the spices having the same gender would not produce anything else while the different gender spices will produce something during the evolution. By $S(\Lambda,G,\Phi),$ we denote the set of all probability distributions of the set $\Omega$. Every element of the simplex $S(\Lambda,G,\Phi)$ can be considered as a state of the population.

We define the evolution of the bisexual population system as a QSO $V:S(\Lambda,G,\Phi)\rightarrow S(\Lambda,G,\Phi)$: during the evolution, every state $\lambda \in S(\Lambda,G,\Phi)$ of the population system goes to another state $\lambda'=V\lambda$ as follows:
\begin{equation}\label{gkso1}
\lambda'(\sigma)=\sum\limits_{\sigma', \sigma'' \in \Omega} p_{\sigma'\sigma'',\sigma}
\lambda(\sigma')\lambda(\sigma''), \quad \ \forall \ \sigma\in\Omega
\end{equation}
where, $p_{\sigma'\sigma'',\sigma}$ is heredity coefficients such that
\begin{equation*}\label{omegaconst}
p_{\sigma'\sigma'',\sigma}=p_{\sigma''\sigma',\sigma}\geq 0
\end{equation*}
\begin{equation*}
\sum\limits_{\sigma \in \Omega(\sigma',\sigma'')} p_{\sigma'\sigma'',\sigma} = \left\{\begin{array}{lll}
0, \ \  \mbox{if} \ \ \Omega(\sigma',\sigma'')= \emptyset \\
c, \ \  \mbox{if} \ \ \Omega(\sigma',\sigma'') \neq \emptyset\\
\end{array}\right.
\end{equation*}
\begin{equation*}
\sum\limits_{\sigma \notin \Omega(\sigma',\sigma'')} p_{\sigma'\sigma'',\sigma}=0.
\end{equation*}
Here $c$ is a given positive constant.

As we mentioned above that depending on possible scenarios for the evolution of the population system, the sex ration might or might not be $\mathbf{1:1}$ (see \cite{Hamil}).

Therefore, in general, we suppose that the population system satisfies the following {\bf $\mathbf{p:q}$ law:} during the evolution, the ratio of females and males remains the same as $\mathbf{p:q}$.
In other words, the population system has an equilibrium $\mathbf{p:q}$ sex ratio.

We are going to provide some models in which one may observe a $\mathbf{p:q}$ law.

We say that an evolution operator \eqref{gkso1} has {\textit{a $\mathbf{p:q}$ property}} if for any $\sigma',\sigma'' \in \Omega$ and $\sigma_f \in \Omega_f,\sigma_m \in \Omega_m$ such that $\sigma_f\mid_{\Lambda\setminus\Lambda_0} = \sigma_m\mid_{\Lambda\setminus\Lambda_0}$ one has
\begin{equation}\label{pq1}
p_{\sigma'\sigma'',\sigma_f}:p_{\sigma'\sigma'',\sigma_m}=p:q,
\end{equation}
\begin{equation}\label{2pq}
\sum\limits_{\sigma \in \Omega} p_{\sigma'\sigma'',\sigma} =
\left\{\begin{array}{l}
0 \ \ \ \ \ \mbox{if} \quad \Omega(\sigma',\sigma'')= \emptyset \\
\frac{1}{2pq} \ \  \mbox{if} \quad \Omega(\sigma',\sigma'') \neq \emptyset,
\end{array}\right.
\end{equation}
where, $0<p,q<1$ and $p+q=1.$

Biological meaning of the condition \eqref{pq1} in the $\mathbf{1:1}$ law case is that the possibilities having the same genotype of girls and boys are equal.

An evolution operator \eqref{gkso1} having a $\mathbf{p:q}$ property is called {\textit{a $\mathbf{p:q}$ operator}}.

Let us define \textit{a hyper-simplex} as follows:
$$S^{*}(\Lambda,G,\Phi)=\{\lambda \in S(\Lambda,G,\Phi): \lambda(\Omega_f)=p,\  \lambda(\Omega_m)=q\}$$

One can easily check that every $\mathbf{p:q}$ operator has the following canonical form
\begin{equation}\label{pq5}
\lambda'(\sigma)=2\sum\limits_{\sigma' \in \Omega_f  \atop \sigma'' \in \Omega_m}  p_{\sigma'\sigma'',\sigma} \lambda(\sigma')\lambda(\sigma'') , \quad \ \forall \ \sigma\in\Omega.
\end{equation}
Indeed, since $\Omega=\Omega_f\cup\Omega_m$ and $\Omega_f\cap\Omega_m=\emptyset,$ we have that
\begin{eqnarray*}
\lambda'(\sigma)&=&\sum\limits_{\sigma',\sigma'' \in \Omega_f \cup \Omega_m} p_{\sigma'\sigma'',\sigma}\lambda(\sigma')\lambda(\sigma'')\\
&=&\sum\limits_{\sigma'\in \Omega_f \atop \sigma'' \in \Omega_m}p_{\sigma'\sigma'',\sigma}\lambda(\sigma')\lambda(\sigma'')+\\
&& \quad \ \ \quad \ \ \quad  +\sum\limits_{\sigma' \in \Omega_m \atop \sigma''\in \Omega_f}p_{\sigma'\sigma'',\sigma}\lambda(\sigma')\lambda(\sigma'')\\
&=&2\sum\limits_{\sigma' \in \Omega_f  \atop \sigma'' \in\Omega_m}p_{\sigma'\sigma'',\sigma}\lambda(\sigma')\lambda(\sigma'').
\end{eqnarray*}

One may check that {\emph{any $\mathbf{p:q}$ operator satisfies a $\mathbf{p:q}$ law.}}

Indeed, we want to show that the hyper-simplex $S^{*}(\Lambda,G,\Phi)$ is invariant under the $\mathbf{p:q}$ operator, i.e., $V\left(S^{*}(\Lambda,G,\Phi)\right)\subset S^{*}(\Lambda,G,\Phi).$

Let us first show that if $\lambda\in S^{*}(\Lambda,G,\Phi)$ then $V\lambda=\lambda'\in S(\Lambda,G,\Phi).$ In fact, it follows from \eqref{2pq} and \eqref{pq5} that
\begin{eqnarray*}
\sum\limits_{\sigma \in \Omega} \lambda'(\sigma)&=& 2 \sum\limits_{\sigma' \in \Omega_f  \atop \sigma'' \in \Omega_m}\bigg(\sum\limits_{\sigma \in \Omega}
p_{\sigma'\sigma'',\sigma}\bigg) \lambda(\sigma')\lambda(\sigma'')\\
&=& \frac{1}{pq} \lambda(\Omega_f)\cdot\lambda(\Omega_m)=1.
\end{eqnarray*}
Now, we want to show that $\lambda'\in S^{*}(\Lambda,G,\Phi).$ Indeed, for any $\sigma_f \in \Omega_f,\sigma_m \in \Omega_m$ such that $\sigma_f|_{\Lambda_1} = \sigma_m|_{\Lambda_1}$, it follows from \eqref{pq1} and \eqref{pq5} that
$\lambda'(\sigma_f):\lambda'(\sigma_m)=p:q.$ This yields that $\lambda'(\Omega_f)=p$ and $\lambda'(\Omega_m)=q.$

In the sequel, we suppose that the $\mathbf{p:q}$ operator acts on the hyper-simplex $S^{*}(\Lambda,G,\Phi)$, i.e., $V: S^{*}(\Lambda,G,\Phi)\to S^{*}(\Lambda,G,\Phi).$  Moreover, for the sake of simplicity, we shall consider the population system having a $\mathbf{1:1}$ law (or Fisher's principle), i.e., $p=q=\frac{1}{2}.$

\section {Mendelian and Non-Mendelian QSO}

In this section, we are going to provide Mendelian and Non-Mendelian models of the bisexual population system having $\mathbf{1:1}$ law (or Fisher's principle).

Mendel's law of the heredity can be summarized in two laws: {\textit{Law of segregation}} and \textit{Law of independent assortment}. Law of segregation states that when any individual produces gametes, the copies of a gene separate so that each gamete receives only one copy allele. Law of independent assortment states that alleles of different genes assort independently of one another during gamete formation.

Now, we want to describe Mendelian and Non-Mendelian models of the bisexual population system.

Let  $\{\Lambda_i\},$ where $i=0,1,...,n$, be a set of maximal connected subgraphs of the graph $(\Lambda,L)$.
For the two genotypes (parents) $\sigma',\sigma'' \in \Omega,$ we define a set (traits of children) $\Omega_{M}(\sigma',\sigma'')$ as follows: if $\sigma'\mid_{\Lambda_0} = \sigma''\mid_{\Lambda_0}$ then $\Omega_M(\sigma',\sigma'')=\emptyset$ and if $\sigma'\mid_{\Lambda_0} \neq \sigma''\mid_{\Lambda_0}$ then $\Omega_M (\sigma',\sigma'')=\{\sigma \in \Omega: \sigma|_{\Lambda_i}=\sigma'|_{\Lambda_i} \ \mbox{or} \ \sigma|_{\Lambda_i}=\sigma''|_{\Lambda_i} \}$.

\begin{definition} \label{defcoef}
\, A $\mathbf{1:1}$ quadratic stochastic operator \eqref{pq5} is called Mendelian, if $p_{\sigma'\sigma'',\sigma}=0,$
for all $\sigma \notin \Omega_M(\sigma',\sigma'')$. A $\mathbf{1:1}$ quadratic stochastic operator \eqref{pq5} is called Non-Mendelian if $p_{\sigma'\sigma'',\sigma_0}\neq 0$
for some $\sigma_0 \notin \Omega_M(\sigma',\sigma'').$
\end{definition}

One of the natural way to construct Mendelian QSO as follows. Let 
$S^{*}(\Lambda,G,\Phi)=\{\lambda \in S(\Lambda,G,\Phi)\lambda(\Omega_f)=\lambda(\Omega_m)=\frac12\}$ be a hyper-simplex.

Let $\mu_{0}\in S^{*}(\Lambda,G,\Phi)$ be a fixed measure such that
\begin{equation}\label{muf=mum}
\mu_0(\sigma_f)=\mu_0(\sigma_m),
\end{equation}
for any $\sigma_f \in \Omega_f,\sigma_m \in \Omega_m$ in which $\sigma_f|_{\Lambda\setminus\Lambda_0} = \sigma_m|_{\Lambda\setminus\Lambda_0}.$

The heredity coefficients $\{p_{\sigma^\prime\sigma^{\prime\prime},\sigma}\}$ of the Mendelian QSO defined as follows
\begin{eqnarray}\label{coef}
p_{\sigma'\sigma'',\sigma}=\begin{cases}
\frac{2\mu_0(\sigma)}{\mu_0\left(\Omega_M(\sigma',\sigma'')\right)} & \mbox{if} \  \sigma \in \Omega_M(\sigma',\sigma'')\neq\emptyset\\
0 &  \mbox{if} \  \  \sigma\not\in\Omega_M(\sigma',\sigma'')\neq\emptyset\\
0 &  \mbox{if} \ \ \Omega_M(\sigma',\sigma'')=\emptyset\\
\end{cases}
\end{eqnarray}
for any $\sigma',\sigma''\in\Omega$.

Then, one can easily check that
\begin{eqnarray*}
p_{\sigma'\sigma'',\sigma_f}:p_{\sigma'\sigma'',\sigma_m}&=&1:1 \\
\sum\limits_{\sigma \in \Omega} p_{\sigma'\sigma'',\sigma}&=&
\left\{\begin{array}{l}
0 \ \ \ \ \mbox{if} \quad \Omega(\sigma',\sigma'')= \emptyset \\
2 \ \ \ \ \mbox{if} \quad \Omega(\sigma',\sigma'') \neq \emptyset,
\end{array}\right.
\end{eqnarray*}
for any $\sigma',\sigma'' \in \Omega$ and $\sigma_f \in \Omega_f,\sigma_m \in \Omega_m$ such that $\sigma_f|_{\Lambda\setminus\Lambda_0} = \sigma_m|_{\Lambda\setminus\Lambda_0}$.

This means that the $\mathbf{1:1}$ QSO with \eqref{coef} is Mendelian.

It is worth mentioning that Mendelian QSO defined by \eqref{coef} is convenient to study dynamics of some bisexual population systems. In the next section we shall consider some applications of Mendelian QSO.

\section {Mendelian QSO on the simplest graph}

\subsection{A Mendelian Trait.} In this section we are going to discuss the mathematical model of a Mendelian trait. {\emph{A Mendelian trait}} is one that is controlled by a single locus and shows a simple Mendelian inheritance pattern. These traits include PTC (phenylthiocarbamide) tasting, hairline shape, tongue rolling, earlobe attachment, hand clasping etc. (see \cite{OMIN}).

Let $(\Lambda,L)$ be a graph such that $\Lambda_0=\{0\}, \Lambda_1=\{1\}, \Lambda=\Lambda_0\cup\Lambda_1$ and $L=\emptyset$. Suppose that a child receives \emph{dominant} and \emph{recessive} alleles from parents, i.e.,  $\Phi_2=\{A,a\}$.  We then have that $\Omega=\{\sigma_1,\sigma_2,\sigma_3,\sigma_4\},$ $\Omega_f=\{\sigma_1,\sigma_2\},$ and $\Omega_m=\{\sigma_3,\sigma_4\}$ where
$$\sigma_1=(f,A), \quad \sigma_2=(f,a), \quad  \sigma_3=(m,A), \quad  \sigma_4=(m,a).$$
For any $\sigma',\sigma''\in\Omega$ one has that
$$
\Omega_M^{\sigma',\sigma''} =
\begin{cases}
\emptyset & \mbox{if } \sigma'\mid_{\Lambda_0}=\sigma''\mid_{\Lambda_0} \\
\{\sigma',\sigma''\} & \mbox{if } \sigma'\mid_{\Lambda_0}\neq\sigma''\mid_{\Lambda_0}, \sigma'\mid_{\Lambda_1}=\sigma''\mid_{\Lambda_1}\\
\Omega & \mbox{if } \sigma'\mid_{\Lambda_0}\neq\sigma''\mid_{\Lambda_0}, \sigma'\mid_{\Lambda_1}\neq\sigma''\mid_{\Lambda_1}
\end{cases}
$$

Let $S^{*}(\Lambda,G,\Phi_2)= \{\lambda \in S(\Lambda,G,\Phi_2): \lambda(\Omega_f) = \lambda(\Omega_m) = \frac12\}$ be a hyper-simplex. In other words, $\lambda\in S^{*}(\Lambda,G,\Phi_2)$ means that $\lambda(\sigma_i)\geq0,$ $i=\overline{1,4}$ such that
$\lambda(\sigma_1)+\lambda(\sigma_2)=\lambda(\sigma_3)+\lambda(\sigma_4)=\frac{1}{2}.$

Now, we choose a measure $\mu_0\in S^{*}(\Lambda,G,\Phi_2)$ satisfying the condition \eqref{muf=mum}, i.e., $\mu_0(\sigma_i)=\alpha_i,$ $i=\overline{1,4}$ such that $\alpha_1=\alpha_3=\alpha$ and $\alpha_2=\alpha_4=\frac{1}{2}-\alpha$ where $\alpha\in (0,\frac{1}{2}).$

We then may define heredity coefficients $\{p_{\sigma^\prime\sigma^{\prime\prime},\sigma}\}$ of the Mendelian $\mathbf{1:1}$  QSO  by the formula (\ref{coef}).

The evolution of the Mendelian trait of the bisexual population system can be given as follows $\lambda'=V\lambda$
\begin{eqnarray}\label{operator2}
\lambda^\prime(\sigma_1)&=&2\lambda(\sigma_1)\lambda(\sigma_3)+4\alpha\lambda(\sigma_2)\lambda(\sigma_3)+\nonumber\\
&& \quad \quad \quad \quad \quad +4\alpha\lambda(\sigma_1)\lambda(\sigma_4)\nonumber\\
\lambda^\prime(\sigma_2)&=&2(1-2\alpha)\lambda(\sigma_1)\lambda(\sigma_4)+2(1-2\alpha)\lambda(\sigma_2)\lambda(\sigma_3)
\nonumber\\
&& \quad \quad \quad \quad \quad +2\lambda(\sigma_2)\lambda(\sigma_4)\nonumber\\
\lambda^\prime(\sigma_3)&=&2\lambda(\sigma_1)\lambda(\sigma_3)+
4\alpha\lambda(\sigma_1)\lambda(\sigma_4)+\\
&&\quad \quad \quad \quad \quad +4\alpha\lambda(\sigma_2)\lambda(\sigma_3)\nonumber\\
\lambda^\prime(\sigma_4)&=&2(1-2\alpha)\lambda(\sigma_1)\lambda(\sigma_4)+2(1-2\alpha)\lambda(\sigma_2)\lambda(\sigma_3)
\nonumber\\
&&\quad \quad \quad \quad \quad +2\lambda(\sigma_2)\lambda(\sigma_4)\nonumber.
\end{eqnarray}

Let
\begin{eqnarray}
S^{2n-1}=\left\{x\in \mathbb{R}^{2n}:  \sum\limits_{i=1}^{2n}x_i=1, \ x_i\geq0, \ \forall i=\overline{1,2n} \right\}
\end{eqnarray}
\begin{eqnarray}\label{hypersimplex}
{\bar{S}}^{2n-1}=\left\{x\in S^{2n-1}: \sum\limits_{i=1}^{n}x_i=\sum\limits_{i=n+1}^{2n}x_i=\frac{1}{2}\right\}
\end{eqnarray}
be a $(2n-1)-$dimensional simplex and a hyper-simplex, respectively.

If we denote by $\lambda(\sigma_i)=x_i$ for all $i=\overline{1,4}$ then the evolution operator $V_\alpha:\bar{S}^{3}\to \bar{S}^{3}$ defined by \eqref{operator2}
has following form
\begin{equation} \label{operator3}
V_\alpha: \left \{\begin{array}{llll}
x^\prime_1= 2x_1x_3 + 4\alpha x_2x_3 + 4\alpha x_1x_4,\\
x^\prime_2= 2(1-2\alpha) x_1x_4 + 2(1-2\alpha) x_2x_3 + 2 x_2x_4,\\
x^\prime_3= 2x_1x_3 + 4\alpha x_2x_3 + 4\alpha x_1x_4,\\
x^\prime_4= 2(1-2\alpha) x_1x_4 + 2(1-2\alpha) x_2x_3 + 2 x_2x_4.
\end{array}\right.
\end{equation}
where $\alpha \in (0,\frac{1}{2})$.

It is clear that for any $x\in\bar{S}^3 $ one holds $(V_\alpha x)_1=(V_\alpha x)_3, \ (V_\alpha x)_2=(V_\alpha x)_4$.

Let us consider the following quadratic function $f_\alpha:[0,\frac12]\to[0,\frac12]$
\begin{equation} \label{function1}
f_\alpha(x)=2(1-4\alpha)x^2+4\alpha x
\end{equation}

One can easily check the following assertions:
\begin{itemize}
  \item[(i) ] The function \eqref{function1} has two fixed points $Fix(f_\alpha)=\left\{0,\frac{1}{2}\right\};$
  \item[(ii) ] If $\alpha=\frac{1}{4}$ then $f_\alpha(x)=x;$
  \item[(iii) ] If $0<\alpha<\frac{1}{4}$ then the trajectory of the function \eqref{function1} converges to the fixed point $x_{*}=0$ for any initial point $x_0\in(0,\frac{1}{2})$. If $\frac{1}{4}<\alpha<\frac{1}{2}$ then the trajectory of the function \eqref{function1} converges to the fixed point $x^{*}=\frac{1}{2}$ for any initial point $x_0\in(0,\frac{1}{2})$.
\end{itemize}

By means of the properties of the quadratic function \eqref{function1} one can get the following properties of the operator $V_\alpha:\bar{S}^{3}\to \bar{S}^{3}$ given by \eqref{operator3}:
\begin{itemize}
  \item[(i) ] The operator \eqref{operator3} has two fixed points  $Fix(V_\alpha)=\left\{\left(0,\frac{1}{2},0,\frac{1}{2}\right),\left(\frac{1}{2},0,\frac{1}{2},0\right) \right\}$;
  \item[(ii) ] If $\alpha=\frac{1}{4}$ then the operator \eqref{operator3} is an identity operator;
  \item[(iii) ] If $0<\alpha<\frac{1}{4}$ then the trajectory of the operator \eqref{operator3} converges to the fixed point $x_{*}=\left(0,\frac{1}{2},0,\frac{1}{2}\right)$ for any initial point $x^0\in\bar{S}^{3}$. If $\frac{1}{4}<\alpha<\frac{1}{2}$ then the trajectory of the operator \eqref{operator3} converges to the fixed point $x^{*}=\left(\frac{1}{2},0,\frac{1}{2},0\right)$ for any initial point $x^0\in\bar{S}^{3}$.
\end{itemize}

Based on this study, we may conclude that,\emph{ in the bisexual population system, the Mendelian trait which is dominating in numbers at the initial state would be spread towards the system in the future.}

\subsection {A Mendelian Inheritance with Multiple Alleles} In this section, we are going to consider a Mendelian inheritance with multiple alleles.

Let $(\Lambda,L)$ be a graph such that $\Lambda_0=\{0\}, \Lambda_1=\{1\}, \Lambda=\Lambda_0\cup\Lambda_1$ and $L=\emptyset$. Suppose that a child receives 4 types of alleles from parents, i.e.,  $ \Phi_4=\{1,2,3,4\}$.

We then have that $$\Omega=\{\sigma_{f,1},\sigma_{f,2},\sigma_{f,3},\sigma_{f,4},\sigma_{m,1},\sigma_{m,2},\sigma_{m,3},\sigma_{m,4}\},$$ 
$$\Omega_f=\{\sigma_{f,1},\sigma_{f,2},\sigma_{f,3},\sigma_{f,4}\},$$
$$\Omega_m=\{\sigma_{m,1},\sigma_{m,2},\sigma_{m,3},\sigma_{m,4}\},$$
where
$$\sigma_{f,1}=(f,1), \ \sigma_{f,2}=(f,2), \ \sigma_{f,3}=(f,3), \ \sigma_{f,4}=(f,4),$$
$$\sigma_{m,1}=(m,1), \ \sigma_{m,2}=(m,2),  \ \sigma_{m,3}=(m,3), \ \sigma_{m,4}=(m,4).$$

For any $\sigma',\sigma''\in\Omega$ one has that
$$
\Omega_M^{\sigma',\sigma''} =
\begin{cases}
\emptyset, & \mbox{if } \sigma'\mid_{\Lambda_0}=\sigma''\mid_{\Lambda_0} \\
\Omega_{\sigma}, & \mbox{if } \sigma'\mid_{\Lambda_0}\neq\sigma''\mid_{\Lambda_0}, \sigma'\mid_{\Lambda_1}=\sigma''\mid_{\Lambda_1}\\
\Omega_{\sigma,\delta}, & \mbox{if } \sigma'\mid_{\Lambda_0}\neq\sigma''\mid_{\Lambda_0}, \sigma'\mid_{\Lambda_1}\neq\sigma''\mid_{\Lambda_1}.
\end{cases}
$$
where, $\Omega_{\sigma}=\{\sigma',\sigma''\}$, $\Omega_{\sigma,\delta}=\{\sigma',\sigma'',\delta',\delta''\}$, and
\begin{eqnarray*}
\delta'\mid_{\Lambda_0}&=&\sigma'\mid_{\Lambda_0}, \  \delta'\mid_{\Lambda_1}=\sigma''\mid_{\Lambda_1}\\
\delta''\mid_{\Lambda_0}&=&\sigma''\mid_{\Lambda_0}, \  \delta'\mid_{\Lambda_1}=\sigma'\mid_{\Lambda_1}.
\end{eqnarray*}

Let $S^{*}(\Lambda,G,\Phi_4)=\{\lambda \in S(\Lambda,G,\Phi_4): \lambda(\Omega_f)=\lambda(\Omega_m)=\frac12\}$ be a hyper-simplex. In other words, $\lambda\in S^{*}(\Lambda,G,\Phi_4)$ means that $\lambda(\sigma_{f,i}),\lambda(\sigma_{m,i})\geq0,$ $i=\overline{1,4}$ such that
$$\sum\limits_{i=1}^4\lambda(\sigma_{f,i})=\sum\limits_{i=1}^4\lambda(\sigma_{m,i})=\frac{1}{2}.$$

Now, we choose a measure $\mu_0\in S^{*}(\Lambda,G,\Phi_4)$ satisfying the condition \eqref{muf=mum}, i.e., $\mu_0(\sigma_{f,i})=\mu_0(\sigma_{m,i})=\alpha_i,$ $i=\overline{1,4}$ such that $\alpha_1+\ldots+\alpha_4=\frac{1}{2}$ and $\alpha_i\geq0$ for any $i=\overline{1,4}.$ We then may define heredity coefficients $\{p_{\sigma^\prime\sigma^{\prime\prime},\sigma}\}$ of the Mendelian $\mathbf{1:1}$  QSO  by the formula (\ref{coef}).

The evolution of the Mendelian inheritance with 4 alleles of the bisexual population system can be given as follows $\lambda'=W\lambda$
\begin{eqnarray*}\label{operator5}
\lambda'(\sigma_{f,i})=\sum\limits_{j=1}^4{\frac{2\alpha_i}{\alpha_i+\alpha_j}\left(\lambda(\sigma_{f,i})\lambda(\sigma_{m,j})
+\lambda(\sigma_{f,j})\lambda(\sigma_{m,i})\right)}\\
\lambda'(\sigma_{m,i})=\sum\limits_{j=1}^4{\frac{2\alpha_i}{\alpha_i+\alpha_j}\left(\lambda(\sigma_{f,i})\lambda(\sigma_{m,j})
+\lambda(\sigma_{f,j})\lambda(\sigma_{m,i})\right)}
\end{eqnarray*} where, $i=\overline{1,4}.$

If we denote by $\lambda(\sigma_{f,i})=x_i$ and $\lambda(\sigma_{m,i})=x_{i+4}$ for all $i=\overline{1,4}$ then the evolution operator $W:\bar{S}^7\rightarrow\bar{S}^7$ defined by \eqref{operator5} has following form
\begin{eqnarray*}\label{operator5}
x'_{i}=\sum\limits_{j=1}^4{\frac{2\alpha_i}{\alpha_i+\alpha_j}}\left(x_{i}x_{j+4}
+x_{j}x_{i+4}\right),\quad i=\overline{1,4},\\
x'_{i+4}=\sum\limits_{j=1}^4{\frac{2\alpha_i}{\alpha_i+\alpha_j}}\left(x_{i}x_{j+4}
+x_{j}x_{i+4}\right),\quad i=\overline{1,4}.
\end{eqnarray*}

It is clear that for any $x\in\bar{S}^7 $ one holds $(Wx)_i=(Wx)_{i+4}$ for all $i=\overline{1,4}.$

Therefore, it is enough to study the following operator
\begin{equation}\label{operator7}
{W}_f:\left \{\begin{array}{lllllllllllllll}
x^\prime_1=x_1\bigg(2x_1+\frac{4\alpha_1}{\alpha_1+\alpha_2}x_2+ \frac{4\alpha_1}{\alpha_1+\alpha_3}x_3+
\frac{4\alpha_1}{\alpha_1+\alpha_4}x_4\bigg),\\
x^\prime_2=x_2\bigg(2x_2+\frac{4\alpha_2}{\alpha_1+\alpha_2}x_1+ \frac{4\alpha_2}{\alpha_2+\alpha_3}x_3+
\frac{4\alpha_2}{\alpha_2+\alpha_4}x_4\bigg),\\
x^\prime_3=x_3\bigg(2x_3+\frac{4\alpha_3}{\alpha_1+\alpha_3}x_1+ \frac{4\alpha_3}{\alpha_2+\alpha_3} x_2+
\frac{4\alpha_3}{\alpha_3+\alpha_4}x_4\bigg),\\
x^\prime_4=x_4\bigg(2x_4+\frac{4\alpha_4}{\alpha_1+\alpha_4}x_1 + \frac{4\alpha_4}{\alpha_2+\alpha_4}x_2+
\frac{4\alpha_4}{\alpha_3+\alpha_4}x_3\bigg).
\end{array}\right.
\end{equation}
on the domain $x_1+\ldots+x_4=\frac{1}{2}$ and $x_i\geq0$ for all $i=\overline{1,4}.$ One may check that the operator \eqref{operator7} has the following form
\begin{equation}\label{operator81}
\left \{\begin{array}{lllllllllllllll}
x^\prime_1=x_1\bigg(1+\frac{2(\alpha_1-\alpha_2)}{\alpha_1+\alpha_2}x_2+ \frac{2(\alpha_1-\alpha_3)}{\alpha_1+\alpha_3}x_3+
\frac{2(\alpha_1-\alpha_4)}{\alpha_1+\alpha_4}x_4\bigg)\\
x^\prime_2=x_2\bigg(1+\frac{2(\alpha_2-\alpha_1)}{\alpha_1+\alpha_2}x_1+ \frac{2(\alpha_2-\alpha_3)}{\alpha_2+\alpha_3}x_3+
\frac{2(\alpha_2-\alpha_4)}{\alpha_2+\alpha_4}x_4\bigg)\\
x^\prime_3=x_3\bigg(1+\frac{2(\alpha_3-\alpha_1}{\alpha_1+\alpha_3}x_1+ \frac{2(\alpha_3-\alpha_2)}{\alpha_2+\alpha_3} x_2+
\frac{2(\alpha_3-\alpha_4)}{\alpha_3+\alpha_4}x_4\bigg)\\
x^\prime_4=x_4\bigg(1+\frac{2(\alpha_4-\alpha_1)}{\alpha_1+\alpha_4}x_1 + \frac{2(\alpha_4-\alpha_2)}{\alpha_2+\alpha_4}x_2+
\frac{2(\alpha_4-\alpha_3)}{\alpha_3+\alpha_4}x_3\bigg).
\end{array}\right.
\end{equation}

Let $y_i=2x_i$ and $y'_i=2x'_i$ for $i=\overline{1,4}.$ Then the operator \eqref{operator81} is quadratic stochastic Volterra operator defined on the simplex $S^{3}:$
\begin{equation}\label{operator8}
y^\prime_i=y_i\bigg(1+\sum\limits_{j=1\atop j\neq i}^4\frac{\alpha_i-\alpha_j}{\alpha_i+\alpha_j}y_j\bigg), \quad i=\overline{1,4}.
\end{equation}

Quadratic stochastic Volterra operators defined on the finite dimensional simplex were studied very well in \cite{GN}-\cite{RSh}. Based on that study, we may get the same conclusion which was done for the Mendelian trait, i.e., \emph{ in the bisexual population system, the dominant Mendelian trait of multiple alleles would be spread towards the system in the future.}

\section {Construction of Non-Mendelian QSO on a finite graph}

In this section, we are aiming to study a special class of non-Mendelian QSO.
This study will be applied to the transmission of Human blood group systems and Rhesus factors (see Appendix).

Let  $\{\Lambda_i\},$ where $i=0,1,...,l$, be a set of maximal connected subgraphs of the graph $(\Lambda,L)$ such that $\Lambda_0=\{x^0\}.$ Let $\Phi$ be a finite set of alleles and $\Omega$ be a set of all genotypes. For the two genotypes (parents) $\sigma',\sigma'' \in \Omega,$ we define a set (traits of children) $\Omega(\sigma',\sigma'')$ as follows:
\begin{eqnarray}\label{Non-mendelianinher}
\Omega(\sigma',\sigma'') =
\begin{cases}
\emptyset, & \mbox{if } \sigma'\mid_{\Lambda_0}=\sigma''\mid_{\Lambda_0} \\
\Omega, & \mbox{if } \sigma'\mid_{\Lambda_0}\neq\sigma''\mid_{\Lambda_0} .
\end{cases}
\end{eqnarray}
In other words, during the evolution of the bisexual population system, a child can receive any type of traits (not only parent's traits).

Now, we are going to provide a construction of a non-Mendelian $\mathbf{1:1}$ QSO defined by a non-Mendelian inheritance \eqref{Non-mendelianinher}.

Let $S^{*}(\Lambda,G,\Phi)=\{\lambda \in S(\Lambda,G,\Phi): \lambda(\Omega_f)=\lambda(\Omega_m)=\frac12\}$ be a hyper-simplex.

For given $\sigma', \sigma''$, we choose a positive measure $\mu_{\sigma',\sigma''}\in S^{*}(\Lambda,G,\Phi)$ such that
\begin{equation}\label{munon-mendelian}
\mu_{\sigma',\sigma''}(\sigma_f)=\mu_{\sigma',\sigma''}(\sigma_m),
\end{equation}
for any $\sigma_f \in \Omega_f,\sigma_m \in \Omega_m$ in which $\sigma_f|_{\Lambda\setminus\Lambda_0} = \sigma_m|_{\Lambda\setminus\Lambda_0}.$

We define the heredity coefficients $\{p_{\sigma^\prime\sigma^{\prime\prime},\sigma}\}$ of the non-Mendelian $\mathbf{1:1}$ QSO as follows
\begin{eqnarray}\label{heredity-Nonmedelian}
p_{\sigma'\sigma'',\sigma}=
\begin{cases}
2\mu_{\sigma',\sigma''}(\sigma) & \mbox{ if } \ \Omega(\sigma',\sigma'')=\Omega,\\
0 & \mbox{ if } \ \Omega(\sigma',\sigma'')=\emptyset.
\end{cases}
\end{eqnarray}

Then the transmission $V:S^{*}(\Lambda,G,\Phi)\rightarrow S^{*}(\Lambda,G,\Phi)$ of Rh blood groups can be written as follows $\lambda'=V\lambda$
\begin{equation}\label{Operator-nonMendelian}
\lambda'(\sigma)=4\sum\limits_{\sigma'\in\Omega_f \atop \sigma''\in \Omega_m} \mu_{\sigma',\sigma''}(\sigma)
\lambda(\sigma')\lambda(\sigma''), 
\end{equation} for all  $\sigma\in\Omega.$

Due to the condition \eqref{munon-mendelian}, one can see that $\lambda'(\sigma_f)=\lambda'(\sigma_m)$ for any $\sigma_f \in \Omega_f,\sigma_m \in \Omega_m$ in which $\sigma_f|_{\Lambda\setminus\Lambda_0} = \sigma_m|_{\Lambda\setminus\Lambda_0}.$ Suppose that $\Omega_f=\left\{\sigma_f^{(1)},...,\sigma_f^{(n)}\right\}$ and $\Omega_m=\left\{\sigma_m^{(1)},...,\sigma_m^{(n)}\right\}$ where $\sigma_f^{(i)}|_{\Lambda\setminus\Lambda_0} = \sigma_m^{(i)}|_{\Lambda\setminus\Lambda_0}$ for all $i=\overline{1,n}.$

Let $p_{ij,k}=2\mu_{\sigma^{(i)}_f, \sigma^{(j)}_m}\left(\sigma^{(k)}_f\right)=2\mu_{\sigma^{(i)}_f, \sigma^{(j)}_m}\left(\sigma^{(k)}_m\right),$ for all $i,j,k=\overline{1,n}.$ Since
\begin{eqnarray*}
\sum\limits_{\sigma\in\Omega}\mu_{\sigma',\sigma''}(\sigma)&=&\sum\limits_{k=1}^n\left(\mu_{\sigma^{(i)}_f, \sigma^{(j)}_m}(\sigma^{(k)}_f)+\mu_{\sigma^{(i)}_f, \sigma^{(j)}_m}(\sigma^{(k)}_m)\right)\\
&=&1
\end{eqnarray*}
we have that $\sum\limits_{k=1}^np_{ij,k}=1$ and $p_{ij,k}\geq0$ for all $i,j,k=\overline{1,n}.$

We denote by $\lambda(\sigma^{(k)}_f)=x_k,$ $\lambda(\sigma^{(k)}_m)=x_{n+k}, \forall k=\overline{1,n}.$
Then the non-Mendelian inheritance evolution \eqref{Operator-nonMendelian} can be written as follows
$$
\begin{cases}\label{Operatornform of onMendelian}
x'_k =2\sum\limits_{i,j=1}^n p_{ij,k}
x_ix_{n+j},\\
x'_{n+k}=2\sum\limits_{i,j=1}^n p_{ij,k}
x_ix_{n+j}.
\end{cases}
$$
for all $k=\overline{1,n}.$ It is clear that $x'_k=x'_{n+k}.$

Let us define $y_k=2x_k$ and $y'_k=2x'_k.$ Then, we may study the following QSO $V:S^{n-1}\to S^{n-1}$ instead of the previous operator
\begin{equation}\label{MainnonMendelian Operator}
y'_k=\sum\limits_{i,j=1}^n p_{ij,k}y_iy_{j},  \ k=\overline{1,n}.
\end{equation}
where, $\sum\limits_{k=1}^np_{ij,k}=1$ and $p_{ij,k}=p_{ji,k}\geq0$ for all $i,j,k=\overline{1,n}.$

The regularity of the QSO \eqref{MainnonMendelian Operator} was studied in \cite{GRAS}. Namely, it was proven that {\emph{if $p_{ij,k}>\frac{1}{2n}$ for all $i,j,k=\overline{1,n}$ then QSO has a unique fixed point and all trajectory converge to this fixed point.}}

The biological interpretation of this result is that, \emph{in the bisexual population system, if initial distributions of all possible traits in the system are greater than the reciprocal of the total number of traits then eventually these distributions become stable and all types of traits will be preserved during the evolution}.

\section*{Acknowledgment}
We thank to anonymous referees for their helpful comments which are leading to the clear presentation of the paper. The Author (U.J.) wishes to thank Faculty of Science, International Islamic University Malaysia, where this paper was written, for offering Postdoctoral fellowships and financial supports.

\appendix
\section{Models of Transmission of ABO and Rh Blood Groups}

In this appendix, we construct QSO which describes inheritance of {\textbf{ABO}}	and Rh blood groups. To the best of our knowledge, this is the first attempt to describe a mathematical model of the inheritance of {\textbf{ABO}}	and Rh blood groups. In our model, we shall not take into consideration any internal or external factors of the society such as migration or immigration. We are considering an ideal society. We could not give any comparison between our outcomes and others because it is the first study on this topic. At the same time, we are not asking any credit to be true our outcomes. Of course, our model is far from being the best model. However, it is the first attempt in this direction and it can be modified and improved for any complex society.   

\subsection {Transmission of Rh Blood Groups.} In this section, we construct QSO which describes inheritance of Rh blood groups.

It is known that there are two types of Rhesus factors: Rh positive and Rh negative, i.e., $\Phi_2=\{+,-\}$. Moreover, a gender of human does not have any influence on the types of the Rhesus factor. Therefore, we can consider a gender of human and Rhesus factors as vertexes of disconnected graph $(\Lambda,L),$ where $\Lambda_0=\{0\}, \Lambda_1=\{1\}, \Lambda=\Lambda_0\cup\Lambda_1$ and $L=\emptyset$. We then have that $\Omega=\{\sigma_{f,+},\sigma_{f,-},\sigma_{m,+},\sigma_{m,-}\},$ $\Omega_f=\{\sigma_{f,+},\sigma_{f,-}\},$ and $\Omega_m=\{\sigma_{m,+},\sigma_{m,-}\}$ where
$$\sigma_{f,+}=(f,+), \quad \sigma_{f,-}=(f,-),$$
$$\sigma_{m,+}=(m,+), \quad  \sigma_{m,-}=(m,-).$$

Now, our main task to calculate heredity coefficients, for example $p_{\sigma_{f,+},\sigma_{m,+}}(\sigma_{m,-}).$ In order to do so, we randomly took around 10,000 parents from Malaysia and investigated the transmission of Rh and \textbf{ABO} blood groups to their children. In our database, $\frac{1}{2}p_{\sigma_{f,+},\sigma_{m,+}}(\sigma_{m,-})$ means a frequency (or von-Mises probability) of having a son with Rh negative from both parents having Rh positive. More precisely, $\frac{1}{2}p_{\sigma_{f,+},\sigma_{m,+}}(\sigma_{m,-})=\frac{N_{{m^{-}}}({f^{+}},{m^{+}})}{N({f^{+}},{m^{+}})},$ where $N({f^{+}},{m^{+}})$ is a total number of children of parents both having Rh positive, and $N_{m^{-}}({f^{+}},{m^{+}})$ is a total number of sons with Rh negative in $N({f^{+}},{m^{+}})$. Similarly, we can calculate other heredity coefficients.     

For the sake of simplicity, we use the following notations
$$\sigma_1=(m,+), \quad \sigma_2=(m,-),$$ 
$$\sigma_3=(f,+), \quad  \sigma_4=(f,-).$$

We define the all heredity coefficients $p_{\sigma^\prime\sigma^{\prime\prime},\sigma}$: 
\begin{equation*}\label{coef5}
\begin{array}{llll}
p_{\sigma_1\sigma_3,\sigma_1}=0.4925 & p_{\sigma_1\sigma_3,\sigma_2}=0.0075 & p_{\sigma_1\sigma_3,\sigma_3}=0.4925 \\
p_{\sigma_1\sigma_3,\sigma_4}=0.0075 & p_{\sigma_1\sigma_4,\sigma_1}=0.3273 & p_{\sigma_1\sigma_4,\sigma_2}=0.1727 \\
p_{\sigma_1\sigma_4,\sigma_3}=0.3273 & p_{\sigma_1\sigma_4,\sigma_4}=0.1727 & p_{\sigma_2\sigma_3,\sigma_1}=0.3230 \\
p_{\sigma_2\sigma_3,\sigma_2}=0.1770 & p_{\sigma_2\sigma_3,\sigma_3}=0.3230 & p_{\sigma_2\sigma_3,\sigma_4}=0.1770 \\
p_{\sigma_2\sigma_4,\sigma_1}=0.05 & p_{\sigma_2\sigma_4,\sigma_2}=0.45 & p_{\sigma_2\sigma_4,\sigma_3}=0.05 \\ 
p_{\sigma_2\sigma_4,\sigma_4}=0.45
\end{array}
\end{equation*}

\begin{remark}
\, The heredity coefficients definitely depend on the region or place where we are collecting statistics. Consequently, the outcomes of distributions of Rh blood groups would be different from place to place. In our statistics, we have randomly chosen 10,000 parents from Malaysia. Moreover, the size of statistics might influence to the outcomes. The outcomes become more accurate and precise as much as large the statistics.        
\end{remark}

Due to \eqref{Operator-nonMendelian}, the transmission of Rh blood groups has the following form
\begin{equation} \label{operator13}
\begin{array}{llll}
x^\prime_1= 1.9699 x_1x_3 + 1.3094 x_1x_4 + 1.2919 x_2x_3 + 0.2 x_2x_4,\\
x^\prime_2= 0.0301 x_1x_3 + 0.6906 x_1x_4 + 0.7081 x_2x_3 + 1.8 x_2x_4,\\
x^\prime_3= 1.9699 x_1x_3 + 1.3094 x_1x_4 + 1.2919 x_2x_3 + 0.2 x_2x_4,\\
x^\prime_4= 0.0301 x_1x_3 + 0.6906 x_1x_4 + 0.7081 x_2x_3 + 1.8 x_2x_4.\\
\end{array}
\end{equation}

It is clear that $x'_1=x'_3$ and $x'_2=x'_4.$ Let $y_1=2x_1$ and $y_2=2x_2.$ As we showed in the previous section, it is enough to study the following QSO $V:S^1\to S^1$ defined on 1-dimensional simplex $S^{1}$
\begin{equation} \label{operatornm3}
V: \left \{\begin{array}{llll}
y^\prime_1= a y_1^2 + 2by_1y_2 + cy_2^2,\\
y^\prime_1= (1-a) y_1^2 + 2(1-b)y_1y_2 + (1-c)y_2^2.\\
\end{array}\right.
\end{equation}
where, $a=0.9849, b=0.6503, c=0.1.$ 

The dynamics of this operator was studied in \cite{Mefa}. Let $\Delta=4(1-a)c+(1-2b)^2$. In our case, one can see that $0<\Delta<4.$
Then the operator \eqref{operatornm3} has a unique attracting fixed point and all its trajectories converge to that fixed point. Moreover, a numerical calculation shows that,
a unique fixed point is approximately equal to $(0.95, 0.05)$ and it is attracting. All trajectory of the operator \eqref{operator13} converge to that fixed point.

Based on this study, we may state that \emph{the transmission of Rh blood groups in Malaysia will be eventually stable and Rh blood groups would be distributed as around 95\% Rh positive and around 5\% Rh negative.}

\subsection {A Transmission of ABO Blood Groups.} In this section, we construct QSO which describes inheritance of \textbf{ABO} blood groups.

It is known that there are four types of \textbf{ABO} blood groups: \textbf{A, B, AB, O}, i.e., $\Phi_4=\{\mathbf{A,B,AB,O}\}$. Moreover, a gender of human does not have any influence on the types of \textbf{ABO} blood groups. Therefore, we can consider a gender of human and \textbf{ABO} blood groups as vertexes of disconnected graph $(\Lambda,L),$ where $\Lambda_0=\{0\}, \Lambda_1=\{1\}, \Lambda=\Lambda_0\cup\Lambda_1$ and $L=\emptyset$. We then have that
$$
\begin{array}{cc}
\Omega=\{\sigma_1,\sigma_2,\sigma_3,\sigma_4,\sigma_5,\sigma_6,\sigma_7,\sigma_8\},\\ \Omega_f=\{\sigma_1,\sigma_2,\sigma_3,\sigma_4\}, \quad \Omega_m=\{\sigma_5,\sigma_6,\sigma_7,\sigma_8\}.
\end{array}
$$
where
$$
\begin{array}{cc}
\sigma_1=(f,{\mathbf{A}}),\ \sigma_2=(f,{\mathbf{B}}), \ \sigma_3=(f,{\mathbf{AB}}), \ \sigma_4=(f,{\mathbf{O}}), \\
\sigma_5=(m,{\mathbf{A}}), \ \sigma_6=(m,{\mathbf{B}}), \ \sigma_7=(m,{\mathbf{AB}}), \ \sigma_8=(m,{\mathbf{O}}).
\end{array}
$$

Now, our main task to calculate heredity coefficients of \textbf{ABO} blood groups, for example $\mu_{\sigma_{f,{\mathbf{A}}},\sigma_{m,{\mathbf{A}}}}(\sigma_{m,{\mathbf{A}}}).$ In our database, $\mu_{\sigma_{f,{\mathbf{A}}},\sigma_{m,{\mathbf{A}}}}(\sigma_{m,{\mathbf{A}}})$ means a frequency (or von-Mises probability) of having a son with {\textbf{A}} blood group  from parents both having {\textbf{A}} blood groups. More precisely, $\mu_{\sigma_{f,{\mathbf{A}}},\sigma_{m,{\mathbf{A}}}}(\sigma_{m,{\mathbf{A}}})=\frac{N_{{m,{{\mathbf{A}}}}}({(f,{{\mathbf{A}}})},{(m,{{\mathbf{A}}})})}{N(({f,{{\mathbf{A}}})},{(m,{{\mathbf{A}}})})},$ where $N(({f,{{\mathbf{A}}})},{(m,{{\mathbf{A}}})})$ is a total number of children of parents both having {\textbf{A}} blood groups, and $N_{{m,{{\mathbf{A}}}}}({(f,{{\mathbf{A}}})},{(m,{{\mathbf{A}}})})$ is a total number of sons with {\textbf{A}} blood group in $N(({f,{{\mathbf{A}}})},{(m,{{\mathbf{A}}})})$. Similarly, we can calculate other heredity coefficients.

Here are the distributions of \textbf{ABO} blood groups in Malaysia: (for the sake of simplicity, we use $\mu_{ij,k}$ instead of $\mu_{\sigma_i\sigma_j}(\sigma_k)$)

\begin{eqnarray*}
\mu_{15,1}=0.4533 & \mu_{16,1}=0.0865 & \mu_{17,1}=0.2218 \\
\mu_{15,2}=0.0063 & \mu_{16,2}=0.2197 & \mu_{17,2}=0.0940 \\
\mu_{15,3}=0.0038 & \mu_{16,3}=0.1661 & \mu_{17,3}=0.1767 \\
\mu_{15,4}=0.0366 & \mu_{16,4}=0.0277 & \mu_{17,4}=0.0075 \\
\mu_{15,5}=0.4533 & \mu_{16,5}=0.0865 & \mu_{17,5}=0.2218 \\
\mu_{15,6}=0.0063 & \mu_{16,6}=0.2197 & \mu_{17,6}=0.0940 \\
\mu_{15,7}=0.0038 & \mu_{16,7}=0.1661 & \mu_{17,7}=0.1767 \\
\mu_{15,8}=0.0366 & \mu_{16,8}=0.0277 & \mu_{17,8}=0.0075 \\
\mu_{25,1}=0.1750 & \mu_{26,1}=0.0060 & \mu_{27,1}=0.0192 \\
\mu_{25,2}=0.1833 & \mu_{26,2}=0.4653 & \mu_{27,2}=0.3846 \\
\mu_{25,3}=0.0983 & \mu_{26,3}=0.0060 & \mu_{27,3}=0.0769 \\
\mu_{25,4}=0.0433 & \mu_{26,4}=0.0227 & \mu_{27,4}=0.0192 \\
\mu_{25,5}=0.1750 & \mu_{26,5}=0.0060 & \mu_{27,5}=0.0192 \\
\mu_{25,6}=0.1833 & \mu_{26,6}=0.4653 & \mu_{27,6}=0.3846 \\
\mu_{25,7}=0.0983 & \mu_{26,7}=0.0060 & \mu_{27,7}=0.0769 \\
\mu_{25,8}=0.0433 & \mu_{26,8}=0.0227 & \mu_{27,8}=0.0192 \\
\mu_{35,1}=0.2525 & \mu_{36,1}=0.0714 & \mu_{37,1}=0.0522 \\
\mu_{35,2}=0.0707 & \mu_{36,2}=0.2662 & \mu_{37,2}=0.0373 \\
\mu_{35,3}=0.1667 & \mu_{36,3}=0.1299 & \mu_{37,3}=0.4030 \\
\mu_{35,4}=0.0101 & \mu_{36,4}=0.0325 & \mu_{37,4}=0.0075 \\
\end{eqnarray*}

\begin{eqnarray*}
\mu_{35,5}=0.2525 & \mu_{36,5}=0.0714 & \mu_{37,5}=0.0522 \\
\mu_{35,6}=0.0707 & \mu_{36,6}=0.2662 & \mu_{37,6}=0.0373 \\
\mu_{35,7}=0.1667 & \mu_{36,7}=0.1299 & \mu_{37,7}=0.4030 \\
\mu_{35,8}=0.0101 & \mu_{36,8}=0.0325 & \mu_{37,8}=0.0075 \\
\mu_{45,1}=0.2730 & \mu_{46,1}=0.0085 & \mu_{47,1}=0.0709 \\
\mu_{45,2}=0.0117 & \mu_{46,2}=0.2888 & \mu_{47,2}=0.0613 \\
\mu_{45,3}=0.0059 & \mu_{46,3}=0.0049 & \mu_{47,3}=0.1142 \\
\mu_{45,4}=0.2094 & \mu_{46,4}=0.1977 & \mu_{47,4}=0.2536 \\
\mu_{45,5}=0.2730 & \mu_{46,5}=0.0085 & \mu_{47,5}=0.0709 \\
\mu_{45,6}=0.0117 & \mu_{46,6}=0.2888 & \mu_{47,6}=0.0613 \\
\mu_{45,7}=0.0059 & \mu_{46,7}=0.0049 & \mu_{47,7}=0.1142 \\
\mu_{45,8}=0.2094 & \mu_{46,8}=0.1977 & \mu_{47,8}=0.2536 \\
\mu_{18,1}=0.2146 & \mu_{28,1}=0.0066 & \mu_{38,1}=0.1070\\
\mu_{18,2}=0.0022 & \mu_{28,2}=0.2325 & \mu_{38,2}=0.1209\\
\mu_{18,3}=0.0045 & \mu_{28,3}=0.0066 & \mu_{38,3}=0.0651\\
\mu_{18,4}=0.2787 & \mu_{28,4}=0.2542 & \mu_{38,4}=0.2070\\
\mu_{18,5}=0.2146 & \mu_{28,5}=0.0066 & \mu_{38,5}=0.1070\\
\mu_{18,6}=0.0022 & \mu_{28,6}=0.2325 & \mu_{38,6}=0.1209\\
\mu_{18,7}=0.0045 & \mu_{28,7}=0.0066 & \mu_{38,7}=0.0651\\
\mu_{18,8}=0.2787 & \mu_{28,8}=0.2542 & \mu_{38,8}=0.2070\\
\mu_{48,1}=0.0058 & \mu_{48,2}=0.0035 & \mu_{48,3}=0.0034\\
\mu_{48,4}=0.4873 & \mu_{48,5}=0.0058 & \mu_{48,6}=0.0035\\
\mu_{48,7}=0.0034 & \mu_{48,8}=0.4873
\end{eqnarray*}

\begin{remark}
\, The heredity coefficients of \textbf{ABO} blood groups definitely depend on the region or place where we are collecting statistics. Consequently, the outcomes of distributions ofof \textbf{ABO} blood groups would be different from place to place. In our statistics, we have randomly chosen 10,000 parents from Malaysia. Moreover, the size of statistics might influence to the outcomes. The outcomes become more accurate and precise as much as large the statistics.        
\end{remark}

Due to \eqref{Operator-nonMendelian}, the transmission $V:S^{7}\rightarrow S^{7}$ of {\textbf{ABO}} blood groups has the following form

\begin{eqnarray*}
x'_1&=&1.8131 x_1x_5+0.3460 x_1x_6+0.8872 x_1x_7+0.858 x_1x_8\\
&+&0.7 x_2x_5+0.0239 x_2x_6+0.0769 x_2x_7+0.0265 x_2x_8\\
&+&1.0101 x_3x_5+0.2857 x_3x_6+0.209 x_3x_7+0.427 x_3x_8\\
&+&1.092 x_4x_5+0.0339 x_4x_6+0.2837 x_4x_7+0.023 x_4x_8
\end{eqnarray*}
\begin{eqnarray*}
x'_2&=&0.0253 x_1x_5+0.8789 x_1x_6+0.3759 x_1x_7+0.009 x_1x_8\\
&+&0.7333 x_2x_5+1.8612 x_2x_6+1.5385 x_2x_7+0.930 x_2x_8\\
&+&0.2828 x_3x_5+1.0649 x_3x_6+0.1493 x_3x_7+0.483 x_3x_8\\
&+&0.047 x_4x_5+1.1554 x_4x_6+0.2452 x_4x_7+0.014 x_4x_8
\end{eqnarray*}
\begin{eqnarray*}
x'_3&=&0.0152 x_1x_5+0.6644 x_1x_6+0.7068 x_1x_7+0.018 x_1x_8\\
&+&0.3933 x_2x_5+0.0239 x_2x_6+0.3077 x_2x_7+0.026 x_2x_8\\
&+&0.6667 x_3x_5+0.5195 x_3x_6+1.6119 x_3x_7+0.260 x_3x_8\\
&+&0.0235 x_4x_5+0.0198 x_4x_6+0.4567 x_4x_7+0.013 x_4x_8
\end{eqnarray*}
\begin{eqnarray*}
x'_4&=&0.147 x_1x_5+0.111 x_1x_6+0.0301 x_1x_7+1.115 x_1x_8\\
&+&0.1733 x_2x_5+0.0909 x_2x_6+0.0769 x_2x_7+1.017 x_2x_8\\
&+&0.0404 x_3x_5+0.1299 x_3x_6+0.0299 x_3x_7+0.828 x_3x_8\\
&+&0.8376 x_4x_5+0.791 x_4x_6+1.0144 x_4x_7+1.949 x_4x_8
\end{eqnarray*}
\begin{eqnarray*}
x'_5&=&1.8131 x_1x_5+0.3460 x_1x_6+0.8872 x_1x_7+0.858 x_1x_8\\
&+&0.7 x_2x_5+0.0239 x_2x_6+0.0769 x_2x_7+0.026 x_2x_8\\
&+&1.0101 x_3x_5+0.2857 x_3x_6+0.209 x_3x_7+0.427 x_3x_8\\
&+&1.092 x_4x_5+0.0339 x_4x_6+0.2837 x_4x_7+0.023 x_4x_8
\end{eqnarray*}
\begin{eqnarray*}
x'_6&=&0.0253 x_1x_5+0.8789 x_1x_6+0.3759 x_1x_7+0.009 x_1x_8\\
&+&0.7333 x_2x_5+1.8612 x_2x_6+1.5385 x_2x_7+0.930 x_2x_8\\
&+&0.2828 x_3x_5+1.0649 x_3x_6+0.1493 x_3x_7+0.483 x_3x_8\\
&+&0.047 x_4x_5+1.1554 x_4x_6+0.2452 x_4x_7+0.014 x_4x_8
\end{eqnarray*}
\begin{eqnarray*}
x'_7&=&0.0152 x_1x_5+0.6644 x_1x_6+0.7068 x_1x_7+0.018 x_1x_8\\
&+&0.3933 x_2x_5+0.0239 x_2x_6+0.3077 x_2x_7+0.027 x_2x_8\\
&+&0.6667 x_3x_5+0.5195 x_3x_6+1.6119 x_3x_7+0.261 x_3x_8\\
&+&0.0235 x_4x_5+0.0198 x_4x_6+0.4567 x_4x_7+0.014 x_4x_8
\end{eqnarray*}
\begin{eqnarray*}
x'_8&=&0.1465 x_1x_5+0.1107 x_1x_6+0.0301 x_1x_7+1.115 x_1x_8\\
&+&0.1733 x_2x_5+0.0909 x_2x_6+0.0769 x_2x_7+1.017 x_2x_8\\
&+&0.0404 x_3x_5+0.1299 x_3x_6+0.0299 x_3x_7+0.828 x_3x_8\\
&+&0.8376 x_4x_5+0.791 x_4x_6+1.0144 x_4x_7+1.949 x_4x_8
\end{eqnarray*}

This operator has a unique fixed point
$$(0.042, 0.258, 0.029, 0.171, 0.042, 0.258, 0.029, 0.171)$$
which is attracting and all its trajectory converge to this fixed point.

Based on this study, we can state that \emph{the transmission of {\textbf{ABO}} blood groups in Malaysia will be eventually stable and {\textbf{ABO}} blood groups would be distributed as follows: around 8\% from {\textbf{A}}, around 52\% from {\textbf{B}}, around 6\% from {\textbf{AB}}, and around 34\% from {\textbf{O}}}.

\end{document}